\documentclass{article}
\usepackage{amsmath}
\usepackage{amsfonts}
\usepackage{amssymb}

\def\N{\mathbb{N}}
\def\R{\mathbb{R}}
\def\d{\delta}
\def\g{\gamma}
\def\n{\noindent}
\def\p{\partial}

\def\C{\mathcal{C}}

\def\w#1{{\widetilde{#1}}}

\newtheorem{lem}{Lemma}
\newtheorem{prop}{Proposition}
\newtheorem{thm}{Theorem}

\newcommand{\be}{\begin{equation}}
\newcommand{\ee}{\end{equation}}
\newcommand{\ba}{\begin{array}}
\newcommand{\ea}{\end{array}}

\begin{document}

\title{Order of tangency between manifolds}

\author{Wojciech Domitrz, Piotr Mormul and Piotr Pragacz}

\date{}

\maketitle

\begin{abstract}
We study the order of tangency between two manifolds of same dimension
and give that notion three quite different geometric interpretations.
Related aspects of the order of tangency, e.\,g., regular separation
exponents, are also discussed.

\vskip3mm
\n{\it Keywords and Phrases.} Order of tangency. Order of contact. 
Taylor polynomial. Higher jets. Tower of Grassmannians. Regular 
separation exponent. \L ojasiewicz exponent.
\vskip2mm
\n 2010 {\it Mathematics Subject Classification.} Primary 14C17, 14M15,
14N10, 14N15, 14H99, 14P10, 14P20, 32B20, 32C07.

\end{abstract}
\section{Introduction}
In the present paper we discuss the order of tangency (or that of contact) between
manifolds and its relation to enumerative geometry started with classical Schubert
calculus.

Two plane curves, both sufficiently smooth and nonsingular at a point $x^0$, 
are said to have a contact of order at least $k$ at $x^0$ if, in properly chosen 
regular parametrizations, those two curves have identical Taylor polynomials of degree $k$ 
about the respective preimages of $x^0$.\footnote{\,Some authors prefer to use at this 
place Taylor polynomials of degree $k-1$ instead, see for instance \cite{CK}.}

Alternatively, those curves have such contact when their minimal regular separation 
exponent at $x^0$, cf. \cite{L}, is strictly bigger than $k$ or is not defined.

Formulas enumerating contacts have been widely investigated. For example in \cite{CK}
the authors derive a formula for the number of contacts of order $n$ between members
of specified $(n-1)$-parameter family of plane curves and a generic plane curve of
a sufficiently high degree.

Contact problems of this sort have been of both old and new interests, particularly 
in the light of Hilbert's 15th problem to make rigorous the classical calculations of 
enumerative geometry, especially those undertaken by Schubert \cite{S}. The situation 
regarding ordinary (i.\,e., second--order) contacts between families of varieties is 
now well understood thanks in large measure to the contact formula of Fulton, Kleiman 
and MacPherson \cite{FKM}. The above mentioned formula in \cite{CK} generalises that 
given by Schubert in \cite{S1} for the number of triple contacts between a given plane 
curve and a specified 2-parameter family of curves. Schubert made his computations 
through the use of what has come to be known as ``Schubert triangles''. This theory 
has been made completely rigorous by Roberts and Speiser, see, e.\,g., \cite{RS1}, 
and independently by Collino and Fulton \cite{CF}.
\vskip.5mm
Apart from contact formulas, an important role is played by the ``order of tangency''.
Let us discuss this notion for Thom polynomials. Among important properties of Thom
polynomials we record their positivity closely related to Schubert calculus (see,
e.\,g., \cite{MPW1} and \cite{P} for a survey).
Namely, the order of tangency allows one to define the jets of Lagrangian submanifolds.
The space of these jets is a fibration over the Lagrangian Grassmannian and leads to
a positive decomposition of the Lagrangian Thom polynomial in the basis of Lagrangian
Schubert cycles.
\vskip1.2mm
In this paper, we give three approaches to the order of tangency. The first one
(in section \ref{un}) is by the Taylor approximations of local parametrizations 
of manifolds. The second one (a mini-max procedure in section \ref{DEUX}) makes use
of curves sitting in the relevant manifolds. The third approach (in section \ref{GRA}) 
is by Grassmann bundles. We show that these three approaches are equivalent. 
We basically work with manifolds over the reals (of various classes of smoothness), 
but the results carry over -- in the holomorphic category -- to complex manifolds. 
\vskip1.2mm
In the last two sections, we discuss some issues related to the ``closeness'' of pairs 
of geometric objects: branches of algebraic sets and relations with contact geometry. 
In fact, in section \ref{quatre} discussed are the regular separation exponents of 
pairs of semialgebraic sets, sometimes called \L ojasiewicz exponents (not to be 
mixed with the by now classical exponents in the renowned \L ojasiewicz inequality 
for analytic functions). Then, in section \ref{con} we report on an unexpected 
application of a modification of tangency order in 3D which yields an elegant 
criterion for a rank-2 distribution on a 3-manifold to be contact.\\
These concluding sections are not less important than the preceding ones. 
They show that the precise measurement of closeness is sometimes more 
demanding -- and giving more -- than merely bounding below tangency orders. 
\vskip1.2mm
In the case of {\it singular\,} varieties different approaches to tangency orders 
lead to different notions. In this respect we refer the reader to \cite{DT} where 
compared were two discrete {\it symplectic\,} invariants of singular curves: 
the {\it Lagrangian tangency order\,} and {\it index of isotropy}. 
\vskip1.5mm
We firstly thank the anonymous referee of \cite{P} for the report which was very 
stimulating for our present studies. Also, we thank Tadeusz Krasi\'nski for informing 
us about the regular separation exponents of pairs of sets, a notion due to \L ojasiewicz. 
Lastly, we thank anonymous referees of the present work for their meticulous 
inspection and advice. 
\section{By Taylor}\label{un}
One situation that is frequently encountered at the crossroads of geometry
and analysis deals with pairs of manifolds which are the graphs of functions
of the same number of variables. Such graphs can intersect, or touch each 
other, at a prescribed point, with various degrees of closeness.
\vskip1mm
Our departing point is a definition of such proximity going precisely in
the spirit of a benchmark reference book \cite{KLV}, p.\,18, although 
not formulated {\it expressis verbis\,} there. 
\vskip2mm
\n{\bf Definition.} Two manifolds $M$ and $\w{M}$ embedded in $\R^m$, both of class 
\,${\rm C}^r$, $r \ge 1$, and the same dimension $p$, intersecting at $x^0 \in M \,\cap \,\w{M}$, 
for $k \le r$, {\it have at $x^0$ the order of tangency at least $k$}, when {\bf there exist} 
a neighbourhood $U \ni u^0$ in \,$\R^p$ and parametrizations\footnote{\,\,the standard topology 
language adopted, among many other sources, in \cite{Mil}} (diffeomorphisms onto the image) 
$$
q \colon \big(U,u^0\big) \to \big(M,\,x^0\big), \ \ \ \ \ \ \tilde{q} 
\colon \big(U,u^0\big) \to \big(\w{M},\,x^0\big)
$$
of class ${\rm C}^r$ such that
\be\label{I}
\Big(\tilde{q} - q\Big)(u) \,= \,{\rm o}\Big(\big|u - u^0\big|^k\Big)
\ee
when $U \ni u \to u^0$. 
\vskip2mm
\n(We underline the existence clause in this definition. Supposing having already 
such a couple of local parametrizations $q$ and $\tilde{q}$, there is {\it an abundance\,} 
of other pairs of \,${\rm C}^r$ parametrizations serving the vicinities of $x^0$ in 
$M$ and $\w{M}$, respectively, and {\it not\,} satisfying the condition \eqref{I}. 
Note also that in this definition the order of tangency is automatically at least 0.)\\
Below in $\bullet\bullet$ in section \ref{DEUX}, and also in section \ref{GRA} 
we restrict ourselves to parametrizations of very specific type -- just the {\it graphs} 
of \,${\rm C}^r$ mappings going from $p$ dimensions to $m - p$ dimensions. This appears 
to be possible while {\it not\,} violating the key condition \eqref{I}. 
\vskip1mm
\n Naturally enough, the notion of the order of tangency not smaller than $\dots$ \,is 
invariant under the local \,${\rm C}^r$ diffeomorphisms of neighbourhoods in \,$\R^m$ 
of the tangency point $x^0$. 
\vskip1.5mm
\n{\it Attention.} In the real \,${\rm C}^\infty$ category it is possible
for the order of tangency to be at least $k$ for all $k \in \N$. In other 
words -- be infinite even though $\{x^0\} = M \,\cap \,\w{M}$. The rest 
of this paper is to be read with this remark in mind.
\vskip2mm
\n As a matter of record, basically the same definition is evoked in 
Proposition on page 4 in \cite{Jensen}. In \cite{Jensen} there is 
also proposed the following reformulation of \eqref{I}.
\begin{prop}\label{or}
The condition \eqref{I} is equivalent to
\be\label{I'}
T_{u^0}^k\big(q\big) \,= \,T_{u^0}^k\big(\tilde{q}\big)\,,
\ee
where $T_{u^0}^k\big(\cdot\big)$ means the Taylor polynomial about $u^0$ of \,degree $k$.
\end{prop}
.\eqref{I} \,$\Rightarrow$ \eqref{I'}.
\begin{eqnarray}\label{I''}
&{\rm o}\Big(\big|u - u^0\big|^k\Big) = \,\tilde{q}(u) - q(u) \,=
\left(\tilde{q}(u) - T_{u^0}^k\big(\tilde{q}\big)(u - u^0)\right)\nonumber \\
&{}+ \,\left(T_{u^0}^k\big(\tilde{q}\big)(u - u^0) - T_{u^0}^k\big(q\big)(u - u^0)\right) +
\left(T_{u^0}^k\big(q\big)(u - u^0) - q(u)\right),
\end{eqnarray}
where the first and last summands on the right hand side are ${\rm o}\Big(\big|u - u^0\big|^k\Big)$ 
by Taylor. Under \eqref{I}, so is the middle summand 
\[
T_{u^0}^k\big(\tilde{q}\big)(u - u^0) - T_{u^0}^k\big(q\big)(u - u^0) \,= \,{\rm o}\Big(\big|u - u^0\big|^k\Big)
\]
and \eqref{I'} follows from the following general result.
\begin{lem}\label{school}
Let $w \in \R[u_1,\,u_2,\dots,\,u_p]$\,, \,${\rm deg}\,w \le k$, $w(u) = {\rm o}\big(|u|^k\big)$
when $u \to 0$ in $\R^p$. Then $w$ is identically zero.
\end{lem}
The proof goes by induction on $k \ge 0$, with an obvious start for $k = 0$. Then, assuming
this for the polynomials of degrees smaller than $k \ge 1$ and taking a polynomial $w$
of degree $k$ as in the wording of the lemma, we can assume without loss of generality
that $w$ is {\it homogeneous\,} of degree $k$ (the terms of lower degrees vanish
altogether by the inductive assumption). Let $u \in \R^p$, $|u| = 1$,be otherwise 
arbitrary. Then
\[
t^k w(u) = w(tu) = {\rm o}\big(|tu|^k\big) = {\rm o}\big(|t|^k\big)\quad
{\rm when\ }t \to 0\,.
\]
Hence $w(u) = 0$ and the vanishing of $w$ follows.
\vskip1.5mm
\n.\eqref{I} $\Leftarrow$ \eqref{I'}.\\
This implication is obvious, because now the middle term on the right hand side of \eqref{I''}
vanishes, so that the right hand side is automatically ${\rm o}\Big(\big|u - u^0\big|^k\Big)$.
\hfill$\Box$
\section{By curves}\label{DEUX}
In the discussion in this section important will be the quantity 
\be\label{s}
s \colon = \,\,\sup\{k \in {\mathbb N} \colon \textrm{the order of tangency}\,\ge k\}\,.
\ee
(Note that an additional restriction here on $k$ is $k \le r$, cf. Definition above.) 
If the class of smoothness $r = \infty$, then, by the very definition, the condition \eqref{I} 
holds for all $k$ if and only if $s = \infty$. 
\vskip1.5mm
\n Is it possible to ascertain something similar in the finite-order-of-tangency case?
\vskip1.5mm
\n With an answer to this question in view, we stick in the present section to
the notation introduced in section \ref{un}, but assume {\it additionally\,} that
\be\label{addit}
s < r\,.
\ee
(Reiterating, the quantity $s$ is defined in \eqref{s} above, and $r$ is the assumed class 
of smoothness of the underlying manifolds, finite or infinite when the category is real. 
When $r = \infty$, the condition \eqref{addit} simply says that $s$ is finite.)
\vskip1.5mm
\n Our second approach uses {\it pairs of curves\,} lying, respectively, 
in $M$ and $\w{M}$.\\We naturally assume that $T_{x^0}M = T_{x^0}\w{M}$. 
Our actual objective is to show that
\begin{thm}\label{deux}
Under \eqref{addit},
\be\label{ddeux}
\underset{v}{\min}\Big(\underset{\g,\,\tilde{\g}}{\max}\big(\max\big\{l \in \{0\} \cup {\mathbb N} 
\,\colon\,|\g(t) - \tilde{\g}(t)| = {\rm o}\big(|t|^l\big)\ {\rm when}\ t \to 0\,\big\}\big)\Big) \,= \,s\,.
\ee
The {\bf minimum} is taken over all \,$0 \ne v \in T_{x^0}M = T_{x^0}\w{M}$. The {\bf outer
maximum} is taken over all pairs of \,\,${\rm C}^r$ curves $\g \subset M$, $\tilde{\g} \subset \w{M}$ 
such that $\g(0) = x^0 = \tilde{\g}(0)$, and -- both non-zero! -- velocities $\dot{\g}(0)$, 
$\dot{\tilde{\g}}(0)$ are both parallel to $v$. 
\end{thm}
{\it Attention}. In this theorem the assumption \eqref{addit} is essential; 
our proof would not work in the situation $s = r$.
\vskip2mm
\n Proof of Theorem \ref{deux}. It is quick to show that the integer on the left hand side 
of equality \eqref{ddeux} is at least $s$. Indeed, for every fixed vector $v$ as above, 
$v = dq(u^0)u$ (without loss of generality, $u$ is like in the proof of 
Lemma \ref{school}), one can take $\d(t) = q(u^0 + tu)$ and 
\,$\tilde{\d}(t) = \tilde{q}\big(u^0 + tu\big)$. Then
\[
|\delta(t) - \tilde{\delta}(t)| \,= \,{\rm o}\big(|tu|^s\big) \,= \,{\rm o}\big(|t|^s\big)
\]
and so 
\[
\underset{\g,\,\tilde{\g}}{\max}\big(\max\big\{\,l \colon |\g(t) - \tilde{\g}(t)|
= {\rm o}\big(|t|^l\big)\ {\rm when}\ t \to 0\big\}\big) \,\ge \,s\,.
\]
In view of the arbitrariness in our choice of $v$, the same remains true after
taking the minimum over all admissible $v$'s which is actually done on the left 
hand side of \eqref{ddeux}.
\vskip2mm
\n$\bullet\bullet$\,To show the opposite non-sharp inequality in \eqref{ddeux} is more
involved. It is precisely in this part that the additional assumption $s \le r - 1$ is
needed. We study the two manifolds in the vicinity of $x^0$ via an appropriate local 
\,${\rm C}^r$ diffeomorphism of the ambient space, after which
\[
\big(M,\,x^0\big) = \Big(\big\{x_{p+1} = x_{p+2} = \dots = x_m = 0\big\},\,0\Big)
\]
and
\[
(\w{M},\,x^0) = \Big(\big\{x_j = F^j(x_1,\,x_2,\dots,\,x_p)\,,\ j = p+1,\,p+2,\dots,\,m\big\},\,0\Big)
\]
for some \,${\rm C}^r$ functions $F^j$. Having the manifolds so neatly 
(graph-like) positioned, we take the most adapted parametrizations 
\[
q(u_1,\,u_2,\dots,\,u_p) = \big(u_1,\,u_2,\dots,\,u_p,\,0,\,0,\dots,\,0\big)\,,
\]
\[
\tilde{q}\big(u_1,\,u_2,\dots,\,u_p\big) = \Big(u_1,\,u_2,\dots,\,u_p,\,F\big(u_1,\,u_2,\dots,\,u_p\big)\Big)\,,
\]
where $F = \big(F^{p+1},\,F^{p+2},\dots,\,F^m\big)$. This --- important --- necessitates 
some extra technical work. Firstly the initial couple of parametrizations satisfying 
\eqref{I} is being straightened simultaneously with manifolds $M$ and $\w{M}$. 
Naturally enough, the resulting parametrizations keep satisfying \eqref{I}, but 
are not yet of the above-desired form. So the parametrizations {\it and\,} manifolds 
are to be additionally slightly upgraded via another local \,${\rm C}^r$ ambient 
diffeomorphism so as (a)\,to keep the simple description of manifolds and 
(b)\,to have the eventual parametrizations adapted as desired above. 
\vskip2mm
\n Given the definition \eqref{s} of $s$, there hold
\[
T_{u^0}^s(q) = T_{u^0}^s\big(\tilde{q}\big)\qquad{\rm and}\qquad T_{u^0}^{s+1}(q)
\ne T_{u^0}^{s+1}\big(\tilde{q}\big)\,,
\]
that is,
\[
T_{u^0}^s(F) = 0\qquad {\rm and}\qquad T_{u^0}^{s+1}(F) \ne 0\,.
\]
It follows that there exist an integer $j \in \{p+1,\,p+2,\dots,\,m\}$
and a vector ${\bf w} \in \R^p$ such that
\be\label{cru}
T_{u^0}^s(F^j\big)({\bf w}) \,= \,0\qquad{\rm and}\qquad T_{u^0}^{s+1}\big(F^j\big)({\bf w}) \ne 0\,.
\ee
Now let $u$ and $\tilde{u}$ be two \,${\rm C}^r$ curves in $\R^p$ passing at $t = 0$ through $u^0$ 
and such that the vectors $\dot{u}(0)$ and $\dot{\tilde{u}}(0)$ are both non-zero and parallel
to {\bf w}. These curves {\it in parameters\,} give rise to \,${\rm C}^r$ curves $\d(t) = q(u(t))$ 
and $\tilde{\d}(t) = \tilde{q}(\tilde{u}(t))$ {\it in the manifolds}, both having at $t = 0$
non-zero speeds parallel to the vector ${\bf v} \colon = dq(u^0){\bf w} = d\tilde{q}(u^0){\bf w}$.
We will now estimate from above (by $s$) the left hand side of the equality \eqref{ddeux}
using, no wonder, {\bf v}, $\d$, and $\tilde{\d}$:
\begin{multline}\label{crux}
|\d(t) - \tilde{\d}(t)| = \sqrt{|u(t) - \tilde{u}(t)|^2 + |F\big(\tilde{u}(t)\big)|^2}\\
{}\ge |F\big(\tilde{u}(t)\big)| \ge |F^j\big(\tilde{u}(t)\big)| \ne {\rm o}\big(|t|^{s+1}\big)\,,
\end{multline}
where the last inequality necessitates an explanation.
In fact, by \eqref{cru} and for every $c \ne 0$
\[
T_{u^0}^{s+1}\big(F^j\big)(tc{\bf w}) = (ct)^{s+1}T_{u^0}^{s+1}\big(F^j\big)({\bf w})
\ne {\rm o}\big(|t|^{s+1}\big)\quad{\rm when}\ \ t \to 0\,.
\]
But $\tilde{u}(t) - \tilde{u}(0) = ct{\bf w} + {\rm o}\big(|t|\big)$ for some non-zero $c$,
hence
\[
T_{u^0}^{s+1}\big(F^j\big)\big(\tilde{u}(t) - \tilde{u}(0)\big) \ne {\rm o}\big(|t|^{s+1}\big)
\quad{\rm when}\ \ t \to 0
\]
as well. Also, just by Peano in the class of smoothness $s+1 \le r$, cf. \eqref{addit},
\[
F^j(u) = T_{u^0}^{s+1}\big(F^j\big)\big(u - u^0\big) + {\rm o}\Big(\big|u - u^0\big|^{s+1}\Big)
\]
when $u \to u^0$ in $\R^p$. Therefore, $F^j\big(\tilde{u}(t)\big)
\ne {\rm o}\big(|t|^{s+1}\big)$ as written in \eqref{crux}.
\vskip1.5mm
Now it is important to note that the pair of curves $\d$ and $\tilde{\d}$ produced by us above 
is completely general in the category \,${\rm C}^r$ for that chosen vector {\bf v}. Hence it 
follows that -- for this precise vector {\bf v}\,! -- the quantity
\[
\underset{\g,\,\tilde{\g}}{\max}\big(\max\big\{\,l \in \{0\} \cup {\mathbb N} \,\colon 
|\g(t) - \tilde{\g}(t)| = {\rm o}\big(|t|^l\big)\ {\rm when}\ t \to 0\big\}\big)
\]
does not exceed $s$. Understandingly, so does the minimum of such quantities
over all $v$'s in $T_{x^0}M = T_{x^0}\w{M}$. Theorem \ref{deux} is now
proved.\hfill$\Box$
\section{By Grassmannians}\label{GRA}
Our third approach is based on the introductory pages of \cite{Jensen} where
a natural tower of consecutive Grassmannians is being attached to every given local
${\rm C}^r$ parametrization $q$ as used by us in the preceding sections. However,
to allow for a recursive definition of tower's members, a more general framework
is needed.
\vskip1.5mm
Namely, to every \,${\rm C}^1$ {\it immersion\,} $H: N \to N'$, $N$ -- an $n$-dimensional
manifold, $N'$ -- an $n'$-dimensional manifold (manifolds not necessarily embedded in 
Euclidean spaces!), we attach the so-called image map \,${\cal G}H \colon N \to G_{n}(N')$ 
of the tangent map $d\,H$: for $s \in N$,
\be\label{gras}
{\cal G}H(s) \,= \,dH(s)(T_sN)\,,
\ee
where $G_{n}(N')$ is the {\it total space} of the Grassmann bundle, with 
base $N'$, of all $n$-planes tangent to $N'$. That is, $G_{n}(N')$ is a new 
manifold, much bigger than $N'$ (whenever $n' > n$), of dimension $n' + n(n' - n)$. 
\vskip2mm
\n We stick in the present section to the notation from section \ref{un} 
and invariably use the pair of parametrizations $q$ and $\tilde q$. 
So we are given the mappings
\[
{\cal G}\,q \,\,\colon \,U \longrightarrow G_p\big(\R^m\big)\,,\qquad
{\cal G}\,{\tilde q} \,\,\colon \,U \longrightarrow G_p\big(\R^m\big)\,.
\]
Upon putting $M^{(0)} = \,\R^m$, \,${\cal G}^{(1)} = \,{\cal G}$, there emerge 
two sequences of recursively defined mappings. \,Namely, for \,$l \ge 1$, 
\[
{\cal G}^{(l)}q \,\,\colon \,U \longrightarrow G_p\big(M^{(l-1)}\big)\,,
\qquad {\cal G}^{(l + 1)}q \,= \,{\cal G}\Big({\cal G}^{(l)}q\Big)
\]
and
\[
{\cal G}^{(l)}{\tilde q} \,\,\colon \,U \longrightarrow
G_p\big(M^{(l-1)}\big)\,,\qquad{\cal G}^{(l + 1)}{\tilde q}
\,= \,{\cal G}\Big({\cal G}^{(l)}{\tilde q}\Big),
\]
where, naturally, $M^{(l)} = \,G_p\big(M^{(l-1)}\big)$. 
\,Now our objective is to show the following.
\begin{thm}\label{trois}
\ ${\rm C}^r$ manifolds $M$ and \,$\w{M}$ have at \,$x^0$ the order of tangency 
at least $k$ ($1 \le k \le r$) iff \,there exist \,\,${\rm C}^r$ parametrizations 
$q$ and $\tilde q$ of the vicinities of \,$x^0$ in, respectively, $M$ and \,$\w{M}$, 
such that 
\be\label{dtrois}
{\cal G}^{(k)}q\,(u^0) \,= \,{\cal G}^{(k)}\tilde{q}\,(u^0)\,.
\ee
\end{thm}
(Observe that, in \eqref{dtrois}, there is clearly encoded that $q(u^0) = x^0 = \tilde{q}(u^0)$.)
\subsection{Proof of Theorem \ref{trois}}\label{pff}
In what follows, of interest for us will be the situations when $H$ in \eqref{gras} above
is locally (and all is local in tangency considerations!) the graph of a \,${\rm C}^1$ mapping 
$h \colon \R^p \supset U \to \R^t$. That is, for $u \in U$, $H(u) = \big(u,\,h(u)\big) 
\in \R^{p + t} = \R^p \times \R^t$. Then \eqref{gras}
assumes by far more precise form
\be\label{nothandy}
{\cal G}H(u) = \Big(u,\,h(u)\,;\;d\big(u,h(u)\big)(u)\Big) =
\Big(u,\,h(u);\ {\rm span}\{\partial_j + h_j(u) \colon j = 1,\,2,\dots,\,p\}\Big)\,
\ee
where the symbol $h_j$ means the partial derivative of a vector mapping $h$ with
respect to the variable $u_j$ ($j = 1,\dots,\,p$), and $\partial_j + h_j(u)$ denotes 
the partial derivative of the vector mapping $(\iota,\,h) \colon U \to \R^p\big(u_1,\dots,\,u_p\big)
\times \R^t$ with respect to $u_j$, where $\iota \colon U \hookrightarrow \R^p$ is the inclusion.
\vskip2mm
\n Now observe that the expression for ${\cal G}H(u)$ on the right hand side of \eqref{nothandy} 
is still not quite useful. Yet there are charts in each newly appearing Grassmannian 
(see, for instance, \cite{Jensen} or p.\,46 in \cite{BN})!
\vskip1.5mm
The chart in a typical fibre $G_p$ over a point in the base $\R^{p + t}$, good
for \eqref{nothandy}, consists of all the entries in the bottommost rows (indexed
by numbers $p + 1,\,p + 2,\,\dots,\,p + t$) in the $(p + t) \times p$ matrices
\[
\left[\begin{array}{ccccccc}
v_1 & | & v_2 & | & \dots & | & v_p
\end{array}\right]
\]
with non-zero upper $p \times p$ minor, {\bf after} multiplying the matrix on the right 
by the inverse of that upper $p \times p$ submatrix. That is to say, taking as the local 
coordinates all the entries in rows $(p + 1)$-st,$\dots$,\,$(p + t)$-th of the matrix
\[
\Big[v_j^{\,\,i}\Big]_{\substack{1 \le i \le p + t\\1 \le j \le p}}
\left(\left[v_j^{\,\,i}\right]_{\substack{1 \le i \le p\\1 \le j \le p}}\right)^{-1}. 
\]
That is, these coordinates are all $t \times p$ entries of the matrix 
\[
\Big[v_j^{\,\,i}\Big]_{\substack{p + 1 \le i \le p + t\\1 \le j \le p}}
\left(\left[v_j^{\,\,i}\right]_{\substack{1 \le i \le p\\1 \le j \le p}}\right)^{-1}. 
\]
In these, extremely useful, glasses the description \eqref{nothandy}
gets stenographed to
\be\label{handy}
{\cal G}H(u) = \Big(u,\,h(u)\,;\ \frac{\partial h}{\partial u}(u)\Big),
\ee
where, under the symbol $\frac{\partial h}{\partial u}(u)$ 
understood are all the entries of this {\it Jacobian\,}\\
$(t \times p)$-matrix written in row and separated by commas. 
This technical simplification is central for a proof that follows.
\vskip2.5mm
After this, basically algebraic, preparation we come back to Theorem \ref{trois}. 
\vskip1mm
\n The order of tangency between $M$ and $\w{M}$ at $x^0$ being at least $k$ 
precisely means (Proposition \ref{or}) the existence of local \,${\rm C}^r$ 
parametrizations $q$ and $\tilde q$ satisfying \eqref{I'}. So we are just 
going to show that \eqref{dtrois} $\Longleftrightarrow$ \eqref{I'}.\\
Moreover, we assume without loss of generality --  much like it has been 
the case in the part $\bullet\bullet$ of the proof of Theorem \ref{deux} -- 
that $M$ and $\w{M}$ are locally the graphs of parametrizations $q$ and $\tilde q$, 
respectively. Which, at the same time, satisfy \eqref{I'}. So \eqref{I'} holds 
for $q(u) = (u, f(u))$, $f \colon U \to \R^{m - p}\big(y_{p+1},\dots,\,y_m\big)$ 
and for 
$\tilde{q}(u) = (u, \tilde{f}(u))$, $\tilde{f} \colon U \to \R^{m - p}\big(y_{p+1},\dots,\,y_m\big)$,
$x^0 = \big(u^0, f(u^0)\big) = \big(u^0, \tilde{f}(u^0)\big)$.
\vskip2.5mm
\n.\eqref{I'} \,$\Rightarrow$ \eqref{dtrois}.

\n We will derive such expressions for ${\cal G}^{(k)}q\,(u)$ and ${\cal G}^{(k)}\tilde{q}\,(u)$,
$u \in U$, that the use of the condition \eqref{I'} will just prompt by itself. 
An added value of this derivation will be the control over the sets of natural local coordinates 
in the Grassmannians in question. (With this information at hand the opposite implication 
$\eqref{I'} \Leftarrow \eqref{dtrois}$ will follow in no time.) \,Our main technical tool 
for the $\Rightarrow$ implication is
\begin{lem}\label{LEMMA}
For $1 \le l \le k$ there exists a local chart on the Grassmannian space 
$G_p\big(M^{(l-1)}\big)$ in which the mapping \,${\cal G}^{(l)}q$ evaluated at $u$
assumes the form
\[
\left(u,\,f(u);\,\binom{l}{1} \times f_{[1]}(u),\ \binom{l}{2} \times f_{[2]}(u),\,\dots,\,
\binom{l}{l} \times f_{[l]}(u)\right),
\]
where $f_{[\nu]}(u)$ is a shorthand notation for the aggregate of {\it all\,} the partials 
of the $\nu$-th order at $u$, of all the components of \,$f$, which are in the number 
$(m - p) \times p^{\,\nu}$, and the symbol $N \times (\ast)$ stands for the $N$ copies 
going in row and separated by commas, of an object $(\ast)$.
\end{lem}
{\it Attention}. In this lemma we purposefully distinguish mixed derivatives 
taken in different orders, simply disregarding the Schwarz symmetricity discovery. 
\vskip2.5mm
\n Proof. $l = 1$. We note that
\[
{\cal G}^{(1)}q\,(u) = \Big(u,\,f(u);\ {\rm span}\{\partial_j + f_j(u) 
\colon j = 1,\,2,\dots,\,p\}\Big)\,,
\]
in the relevant Jensen-Borisenko-Nikolaevskii chart, is nothing but
\[
\big(u,\,f(u);\,f_{[1]}(u)\big) \,= \,\left(u,\,f(u);\,\binom{l}{1} \times f_{[1]}(u)\right).
\]
The beginning of induction is done. 
\vskip2.5mm
\n $l \Rightarrow l + 1$, $l < k$. The mapping ${\cal G}^{(l)}q \colon U \to M^{(l)}$,
evaluated at $u$, is already written down, in appropriate local chart assumed to exist
in $M^{(l)}$, as
\be\label{graduallll}
\Bigg(u, \,f(u), \,\binom{l}{1} \times f_{[1]}(u), \,\binom{l}{2} \times f_{[2]}(u),\dots,\,
\binom{l}{l} \times f_{[l]}(u)\Bigg).
\ee
We work with ${\cal G}^{(l+1)}q = {\cal G}\Big({\cal G}^{(l)}q\Big)$. Now, \eqref{graduallll}
being clearly of the form $H(u) = \big(u,\,h(u)\big)$ in the previously introduced notation,
the mapping $h$ reads
\[
h(u) \,= \,\Big(f(u), \,\binom{l}{1} \times f_{[1]}(u), \,\binom{l}{2} \times f_{[2]}(u),\dots,\,
\binom{l}{l} \times f_{[l]}(u)\Big).
\]
In order to have ${\cal G}H(u)$ written down, in view of \eqref{handy}, one ought to write in row:
$u$, then $h(u)$, and then all the entries of the Jacobian matrix $\frac{\partial h}{\partial u}(u)$,
also written in row and separated by commas. The {\bf latter}, in our shorthand notation,
are computed immediately. Namely
\[
\frac{\partial h}{\partial u}(u) \,= \Big(\binom{l}{0} \times f_{[1]}(u), \,\binom{l}{1} \times f_{[2]}(u),
\,\binom{l}{2} \times f_{[3]}(u),\dots,\,\binom{l}{l} \times f_{[l+1]}(u)\Big).
\]
These entries on the right hand side  are to be juxtaposed with the {\bf former} entries
$\big(u,\,h(u)\big)$. For better readability, we put together the groups of {\it same\,}
partials (a yet another permutation of Grassmann-type coordinates, cf. the wording of the lemma).
In view of the elementary identities $\binom{l}{\nu - 1} + \binom{l}{\nu} = \binom{l+1}{\nu}$,
we get in the outcome
\[
\!\!\!\!\!\!\!\!\!\Big(u, \,f(u), \;\binom{l+1}{1} \times f_{[1]}(u), \;\binom{l+1}{2} \times
f_{[2]}(u),\dots,\;\binom{l+1}{l} \times f_{[l]}(u),\; \binom{l+1}{l+1} \times f_{[l+1]}(u)\Big).
\]
Lemma \ref{LEMMA} is now proved by induction.
\vskip2.5mm
We now take $l = k$ in Lemma \ref{LEMMA} and get, for arbitrary $u \in U$, two similar
visualisations of \,${\cal G}^{(k)}q\,(u)$ and \,${\cal G}^{(k)}\tilde{q}\,(u)$. At that,
the equality \eqref{I'} holds true at $u = u^0$. As a consequence, \eqref{dtrois} follows.
\vskip2mm
\n.\eqref{I'} $\Leftarrow$ \eqref{dtrois}.
\vskip1mm
\n With the information on superpositions of the mappings $\cal G$, gathered
in the course of proving the implication \eqref{I'} $\Rightarrow$ \eqref{dtrois},
this opposite implication is clear. Theorem \ref{trois} is now proved.\hfill$\Box$
\section{Algebraic geometry examples and regular separation exponents}\label{quatre}
In the present section, we shall work with the regular separation exponents 
of pairs of sets, a notion due to \L ojasiewicz \cite{L}. We shall compute these 
exponents in several natural examples, and compare the results with the information 
(when available) about the relevant orders of tangency. These examples deal with 
branches of algebraic sets which often happen to be tangent one to another, 
with various degrees of closeness. 
\vskip2mm
\n{\bf Example 1.} In the work \cite{CCKK} (Figure 2 on page 37 there) 
analyzed is the following algebraic set in $\R^2(x,y)$ 
\be\label{Ken}
\C = \{(x,y) \colon (y - x^2)^2 = x^5\}\,. 
\ee
The two branches of $C$ issuing from the point $(0,0)$,
\[
C_- = \{y = x^2 - x^{5/2}\,,\ x \ge 0\} \ \ \ \hbox{and} \ \ \ C_+ = \{y = x^2 + x^{5/2}\,,\ x \ge 0\}\,,
\]
could be naturally extended to one-dimensional manifolds $D_-$ and $D_+$,
both of class \,${\rm C}^2$ -- the graphs of functions
\[
y_-(x) = x^2 - |x|^{5/2} \ \ \ \ \hbox{and} \ \ \ \ y_+(x) = x^2 + |x|^{5/2}\,,
\]
respectively. The Taylor polynomials of degree 2 about $x = 0$ of $y_-$ and $y_+$
coincide. Hence $D_-$ and $D_+$ have at $(0,0)$ the order of tangency
at least 2 (cf. section \ref{un}), and clearly {\it not\,} at least 3.
\vskip2mm
This example clearly suggests that, in real algebraic geometry,
it would be pertinent to use non-integer measures of closeness.
For instance, for the above sets $y_-(x)$ and $y_+(x)$, we may take
\[
\sup\{\alpha > 0\,\colon\, y_{+}(x) - y_{-}(x) = {\rm o}(|x|^\alpha)\ {\rm when}\ x \to 0\}\,.
\]
This kind of a {\it generalized\,} order of tangency would be 5/2 
in the Colley-Kennedy example.\\
In the local analytic geometry there {\it is\,} a precise name for this notion 
-- the minimal regular separation exponent for two (semialgebraic) branches, say 
$X$ and $Y$, of an algebraic set. That is, here the minimal exponent for $X = C_-$ 
and \,$Y = C_+$. Some authors working after \cite{L} used to call such a quantity 
the {\it \L ojasiewicz exponent} of, in this case, $C_-,\,C_+$ at \,$(0,0)$. 
And denoted it -- when specialized to the present situation -- by 
${\cal L}_{(0,0)}\big(C_-,\,C_+\big)$.\footnote{\,\,However, this 
terminology is not yet definitely settled, as shown in a recent work \cite{KSSz}. 
The authors of the latter speak just descriptively about `the \L ojasiewicz 
exponent for the regular separation of closed semialgebraic sets'.}
\vskip2mm
\n{\bf Remark 1.} (i)\,\,Therefore, in Example 1, the rational number 5/2
is the minimal regular separation exponent of the semialgebraic sets \,$C_{-}$
and \,$C_{+}$ which touch each other at $(0,0)$. 
\vskip1mm
\n(ii)\,\,Example 1 quickly generalizes, by means of the equation $(y - x^N)^2 = x^{2N + 1}$
with $N$ arbitrarily large, to yield a pair of \,${\rm C}^N$\,manifolds having the order 
of tangency at least $N$ and not at least $N + 1$, and having the minimal regular 
separation exponent $N + \frac12$.
\vskip2mm
\n It has not been difficult in Example 1 to discern the pair of branches $C_{-}$ 
and $C_{+}$, initially slightly hidden in a synthetic equation \eqref{Ken}. 
But it can happen considerably worse in this respect. Consider, for instance 
\vskip2mm
\n{\bf Example 2.} (a) an algebraic set in the plane $\R^2(x,y)$ defined by 
a single equation
\be\label{Ar}
(xy)^2 \,= \,\frac14\Big(x^2 + y^2\Big)^3. 
\ee
This set possesses a pair (even more than one such pair) of semialgebraic branches 
touching each other at the point $(0,0)$. Yet it is not so immediate to ascertain 
their minimal regular separation exponent. Only after recognizing in \eqref{Ar} 
the classical {\it quatrefoil\,} $x = \cos(\varphi)\sin(2\varphi)$, 
$y = \sin(\varphi)\sin(2\varphi)$, it becomes quick to compute 
the relevant minimal regular separation exponent equal to 2. 
\vskip1.5mm
\n(b)\,It is even more interestingly with another algebraic set in 2D 
given by the equation 
\be\label{Arr}
\Big(x^2 + y^2 - \frac12 x\Big)^2 = \,\frac14\big(x^2 + y^2\big).
\ee
This set possesses as well a pair of semialgebraic branches $\{y \le 0\}$ and 
$\{y \ge 0\}$ touching each other at the point $(0,0)$. Yet it takes some time 
to find their minimal regular separation exponent. In fact, after discovering 
in \eqref{Arr} the classical {\it cardioid\,} $r = \frac12(1 + \cos\varphi)$, 
that exponent turns out to be -- one more time -- a non-integer (3/2, in the 
occurrence).  
\vskip.8mm
\n(c)\,It is worthy of note that for both explicitly defined algebraic sets \eqref{Ar} 
and \eqref{Arr} one can apply general type {\it bounds above}, for the minimal separation 
exponent, produced in \cite{KSSz}. Yet the estimations got in that way are unrealistically 
(by factors of thousands) high. 
\vskip2.5mm
Returning to the notion of the order of tangency, in the realm of algebraic geometry 
the distinction between that order and minimal regular separation exponent sometimes 
happens to be fairly clear, with both discussed quantities effectively computable. 
An instructive instance of such a situation occurs in \cite{T} (Example 3.5 there). 
\vskip2mm
\n{\bf Example 3.} The author of \cite{T} deals there with a pair of one-dimensional 
algebraic manifolds $N$ and $Z$ in $\R^2(x,y)$ intersecting at $(0,0)$. The manifold 
$N = \{y = 0\}$ is already utmostly simplified, whereas $Z = \{y^d + yx^{d - 1} + x^s = 0\}$ 
depends on two integer parameters $d$ and $s$, $1 < d < s$, $d$ odd. What we are going 
to discuss here is a kind of reworking of Tworzewski's original approach, see also 
{\it Attention\,} below. 
\vskip1.2mm
These manifolds have at $(0,0)$ the order of tangency at least $s - d$, and not at least 
$s - d + 1$, while their minimal regular separation exponent at $(0,0)$ is $s - d + 1$.
\vskip1.5mm
Indeed --- to justify this one tries to present $Z$ as the graph of a function $y = y(x)$.
Clearly, $y(0) = 0$ and a function $y(x)$ could not be divisible by, for instance, $x^{s+1}$. 
So, with no loss of generality,
\[
y(x) \,= \,x^kz(x) - \,x^{s - d + 1}
\]
for certain integer $k \ge 1$ and another function $z(x)$ such that $z(0) \ne 0$.

\n$\bullet$ The possibility $k < s - d + 1$ boils rather quickly down to $k = 1$,
and then\\
to the relation $\big(z - x^{s - d}\big)^d + z = 0$, impossible at $x = 0$, 
for $\big(z(0)\big)^d + z(0) \ne 0$, $d$ being odd.

\n$\bullet\bullet$ So $k \ge s - d + 1$ and now
\[
y(x) \,= \,x^{s - d + 1}z(x) - \,x^{s - d + 1}
\]
for a certain function $z(x)$. Upon substituting this $y(x)$ to the defining
equation of $Z$ and simplifying, $(z - 1)^d\,x^{(d - 1)(s - d)} + z = 0$.
Hence $z(0) = 0$. The Implicit Function Theorem is applicable here
around $(0,0)$, because
\[
\left.\frac{\p}{\p\,z}\Big((z - 1)^dx^{(d-1)(s-d)} + z\Big)\right|_{(0,0)} = \,\,1\,.
\]
One gets a locally unique \,${\rm C}^\infty$ function $z(x)$, $z(0) = 0$, 
hence also a locally unique function $y(x) = x^{s - d + 1}z(x) - \,x^{s - d + 1}$ 
whose graph is $Z$. Because the function $z$ vanishes at $0$, the minuend 
in this expression for $y$ is an `o'of the subtrahend when $x \to 0$. 
So the statements about the order of tangency and minimal regular 
separation exponent follow immediately.
\vskip1.5mm
\n{\it Attention.} When $2|d$, the above-found resolving function $y(x)$
is {\it not\,} the only solution to the defining equation of $Z$. Namely,
the necessary equality $\big(z(0)\big)^d + z(0) = 0$, $z(0) \ne 0$, is
then possible with $z(0) = -1$ and
\[
\left.\frac{\p}{\p\,z}\left(\big(z - x^{s - d}\big)^d + z\right)\right|_{(0,-1)}
= \,\,d(-1)^{d - 1} + 1 \,= \,1 - d \,\ne \,0\,.
\]
Hence the Implicit Function Theorem gives this time a locally unique 
(for that $k = 1$) \,${\rm C}^\infty$ function $\tilde{z}(x)$, $\tilde{z}(0) = -1$.
Then the graph of
\[
\tilde{y}(x) \,= \,x\,\tilde{z}(x) - \,x^{s - d + 1} = \,- x + (\textrm{higher powers of $x$})
\]
is a second branch of $Z$ passing through $(0,0) \in \R^2$, transversal to $N$,
in a stark distinction to the previously found, tangent to $N$, branch.
\vskip2mm
\n{\bf Remark 2.} More generally, one could {\it not\,} hope to get a precise
information regarding the minimal regular separation exponent for the pair of 
manifolds $M,\,\w{M}$ on the sole basis of the assumptions in Theorem \ref{deux}.
That is, basically, under \eqref{addit}. Despite the inequality \eqref{crux},
that exponent need not necessarily be $s + 1$. Following Example 1 earlier 
in this section, one could just take the curves $C_-$ and $C_+$ as the curves
$\d$ and $\tilde{\d}$, respectively, in the proof of Theorem \ref{deux}.
That is, to take $M = C_- = \d$ and $\w{M} = C_+ = \tilde{\d}$. Then, as the reader 
has lately seen, $r = 2$ and the quantity $s$ defined in \eqref{s} is also 2, 
while the minimal regular separation exponent is but $s + \frac12 \,(\,= \,\frac52)$.
\vskip1.5mm
Even restricting oneself to a benchmark setting $\dim M = \dim\w{M}$, an additional 
enormous complication could come from the fact that the intersection $M \cap \w{M}$ 
might be a topologically highly nontrivial set (think about the \,${\rm C}^\infty$ 
category). And it is precisely $M \cap \w{M}$ which enters the definition of regular 
separation exponents for the pair $M,\,\w{M}$.
\section{Relation with contact topology}\label{con}
Unsurprisingly, the notion of order of contact proves useful not only in algebraic
geometry (cf. Introduction), but also in geometry {\it tout court}. One not so obvious
application in the {\it real\,} category deals with the real contact structures in
three dimensions. Our summarising it here follows closely Section 1.6 in \cite{Geiges}.
The author considers there a couple \,$\Sigma \subset M$, where $M$ is a {\it contact\,} 
3-dimensional manifold and \,$\Sigma$ -- a fixed embedded surface in it. 
Contact means $M$ being endowed with a contact structure, say $\xi$, in $TM$.

When one approaches a given point $p \in \Sigma$ by points $q$ staying within $\Sigma$,
a natural question is about the order of smallness of the angle $\angle\big(T_q\Sigma,\,\xi_q\big)$.
If that angle is an `O' of the distance of $q$ to $p$ taken to power $k$ (the distance 
measured in any chosen, and hence every, set of smooth local coordinates about $p$), 
then it is said that $\xi$ has the order of contact at least $k$ with $\Sigma$ at $p$.
(Therefore, what is discussed in this section differs a little from the notion
of closeness of a pair of manifolds investigated in the preceding sections.
Yet the added value is substantial.)

That is, to say that the {\it new} order of contact is at least 1 at a given point $p$ 
is tantamount to saying that $\xi_p = T_p\Sigma$. And it is exactly 0 at $p$ whenever
$\xi_p \ne T_p\Sigma$.
\vskip1mm
\n So it comes as a not small surprise that this elementary notion allows one to
characterise the contact structures as such! Namely, a theorem proved in \cite{Geiges}
asserts that a rank-2 tangent distribution $\xi$ on a 3-dimensional $M$ is contact iff
$\xi$ has the new order of contact {\bf at most} 1 with every surface $\Sigma$ embedded
in $M$, and this at every point of \,$\Sigma$.
\vskip1.5mm
The next natural question in this direction is whether it is possible to similarly
characterise contact structures on $(2n + 1)$-dimensional manifolds, $n \ge 2$.
The author of \cite{Geiges} says nothing in this respect.


{\small \obeylines \parindent0pt
{\ }
Wojciech Domitrz
Faculty of Mathematics and Information Science
Warsaw University of Technology
Koszykowa 75, 00-662 Warszawa, Poland
E-mail: domitrz@mini.pw.edu.pl
{\ }
Piotr Mormul
Institute of Mathematics, University of Warsaw
Banacha 2, 02-097 Warszawa, Poland
E-mail: mormul@mimuw.edu.pl
{\ }
Piotr Pragacz
Institute of Mathematics, Polish Academy of Sciences
\'Sniadeckich 8, 00-656 Warszawa, Poland
E-mail: P.Pragacz@impan.pl
}
\end{document}